\documentclass[twoside]{article}
\usepackage{amsmath,amssymb,amsthm,latexsym,bm,indentfirst,wasysym}
\usepackage{graphicx,epsfig,amscd,mathrsfs,multirow,bm}
\usepackage{cite}
\usepackage{caption}
\usepackage{color}
\usepackage[colorlinks,linkcolor=black,anchorcolor=black,citecolor=blue]{hyperref}
\numberwithin{equation}{section}
\usepackage{titlesec}
\usepackage{tikz}
\titleformat{\section}{\large\bfseries}{}{0pt}{}

\input xy
\xyoption{all}
\input scrload.tex

\allowdisplaybreaks[4]

\catcode`@=11
\long\def\@makefntext#1{\noindent #1}
\newtheorem{theorem}{Theorem}[section]

\newtheorem{definition}[theorem]{Definition}
\newtheorem{remark}[theorem]{Remark}

\newtheorem{example}[theorem]{Example}

\newtheorem{corollary}[theorem]{Corollary}

\newcommand{\edge}{\ar@{-}}

\newcommand{\pf}{\noindent\begin {proof}}
\newcommand{\epf}{\end{proof}}

\newcommand{\Hom}{\mbox{\rm Hom}}

\def\mod{\mathop{\rm mod}\nolimits}

\def\add{\mathop{\rm add}\nolimits}

\def\Hom{\mathop{\rm Hom}\nolimits}

\def\End{\mathop{\rm End}\nolimits}

\def\@evenfoot{}\def\@oddfoot{}

\def\@evenhead{\hbox to\textwidth{\small\rm\thepage \hfill
{\it Xin  Li, Hanpeng  Gao}}} 

\def\@oddhead{\hbox to \textwidth{\small{\it
ICE-closed subcategories and epibricks over one-point extensions
} \hfill\thepage}}   

\catcode`@=12

\def\bc{\begin{center}}
\def\ec{\end{center}}
\def\no{\noindent}
\def\hang{\hangindent\parindent}
\def\textindent#1{\indent\llap{\qquad #1\ \ \enspace}\ignorespaces}
\def\ref{\par\hang\textindent}

\begin{document}


\abovedisplayskip=6pt plus 1pt minus 1pt \belowdisplayskip=6pt
plus 1pt minus 1pt
\thispagestyle{empty} \vspace*{-1.0truecm} \noindent
\vskip 10mm \bc{\Large\bf ICE-closed subcategories and epibricks over one-point extensions   
\footnotetext{\footnotesize
Supported by the National Natural Science Foundation of China (Grant No. 11871071, 12371015, 12301041) and the Science Fundation for Distinguished Young Scholars of Anhui Province (Grant No. 2108085J01).\\
* Corresponding author\\  
E-mail address: lxin909@163.com(Xin LI); hpgao07@163.com(Hanpeng GAO)
} } \ec  

\vskip 5mm
\bc{\bf Xin Li,\ \ \ Hanpeng Gao$^*$}\\  
{\small\it \hspace{-0.95cm}School of Mathematical Sciences, Anhui University, Hefei $230601$, P.R. China
}\ec   

\vskip 1 mm

{\narrower\noindent{\small {\small\bf Abstract}\ \ Let $B$ be the one-point extension algebra of $A$ by an $A$-module $M$. We proved that every ICE-closed subcategory in $\mod A$ can be extended to be some ICE-closed subcategories in $\mod B$.
		In the same way, every epibrick in $\mod A$ can be extended to be some epibricks in $\mod B$.
		The number of ICE-closed subcategories in $\mod B$ and the number of ICE-closed subcategories in $\mod A$ are denoted respectively as $m$, $n$.
		We can conclude the following inequality:
		$$m \geq 2n$$
		This is the analogical in epibricks.
		As an application, we can get some wide $\tau$-tilting modules of $B$ by wide $\tau$-tilting modules of $A$.

\vspace{1mm}\baselineskip 12pt

\no{\small\bf Keywords} \ \ ICE-closed subcategory; epibrick; one-point extension 

\par

\vspace{2mm}

\no{\small\bf MR(2020) Subject Classification\ \ {\rm 16G20}} 

}}

\baselineskip 15pt

\pagestyle{myheadings}
\markboth{\rightline {\scriptsize  Xin LI and Hanpeng GAO}}
         {\leftline{\scriptsize  ICE-closed subcategory,epibrick and one-point extension }}

\section{1. Introduction} 

Several kinds of subcategories have been researched in the representation theory of algebras.
For example, torsion class and torsion-free class are the key points of these subcategories. 
Torsion class is closed under quotients and extensions and can be classified by support $\tau$-tilting modules in \cite{2014tiltingtheory}, which is an important breakthrough in classification of these subcategories.
Similarly, Haruhisa Enomoto's paper given us a uniform way to classify torsion-free class by considering the information on monobricks.

Bricks and semibricks are considered in \cite{1962france},\cite{1976re}. Moreover, the semibrick has been studied from the point of view of $\tau$-tilting theory in \cite{2020semibrick}. In 2021, Haruhisa Enomoto given the definition of monobrick in \cite{2021monobrick}: a set of bricks where every non-zero map between elements of bricks' isomorphism classes is an injection.
The set of simple objects provides an effective approach to investigate torsion-free class and wide subcategories.
Because monobricks are in bijiection with left Schur subcategories, which are same as subcategories closed under kernels, images and extensions.
Without using  $\tau$-tilting theory, it infers several noted consequence on torsion class and wide subcategories via monobricks.

In 2022, the concept of ICE-closed subcategories of module categories have been introduced by Haruhisa Enomoto in \cite{2022ICE-closed}.
The ICE-closed subcategory closed under images,cokernels and extensions correlates closely with torsion class and wide subcategory.
It is worthy mentioning that representative instances of ICE-closed subcategory are torsion class and wide subcategory .
Therefore, the ICE-closed subcategory can be seen as generation of these two classes.
She proved that the number of ICE-closed subcategories does not dictated by the orientation of the quiver and given a clear formula for each Dynkin type.

In this paper, we construct ICE-closed subcategories and epibricks over the one-point extension $B$ of an algebra $A$ by an $A$-module $M_{A}$.
The following is our main results of this article.

\begin{theorem} {\sl
Let $B$ be the one-point extension algebra of $A$ by an $A$-module $M_{A}$ and $\mathcal{T}_{A}$ be an ICE-closed subcategory in $\mod A$.
\begin{enumerate}
	\item[(1)]$\mathcal{T}_{B} := \{ (N_{A},0,0) | N_{A} \in \mathcal{T}_{A}\}$, $\mathcal{T}_{B}$ is an ICE-closed subcategory in $\mod B$.
	\item[(2)]$\mathcal{T}_{B} := \{ (N_{A},k^{n},f) ,(0,k^{n},0)| N_{A} \in \mathcal{T}_{A},n \in \mathbb{N}, f: k^{n}\otimes_k M_{A}\rightarrow N_{A}\}$, $\mathcal{T}_{B}$ is an ICE-closed subcategory in $\mod B$.
\end{enumerate}
}
\end{theorem}
\begin{theorem} {\sl
	Let $B$ be the one-point extension algebra of $A$ by an $A$-module $M_{A}$ and $\mathcal{S}_{A}$ be an epibrick in $\mod A$.
	
	\begin{enumerate}
		\item[(1)]$\mathcal{S}_{B} := \{ (s,0,0) | s \in \mathcal{S}_{A}\}$, $\mathcal{S}_{B}$ is an epibrick in $\mod B$.
		\item[(2)]$\mathcal{S'}_{B} := \{ (s,0,0) ,(0,k,0)| s \in \mathcal{S}_{A}\}$, $\mathcal{S}_{B}$ is an epibrick in $\mod B$.
	\end{enumerate}
}
\end{theorem}

Throughout this paper, all algebras will be basic,  connected, finite dimensional $K$-algebras over an algebraically closed field $K$. Let $A$ be an algebra, $\mod A$ will be the category of finitely generated right $A$-modules and  $\tau$ the Auslander-Reiten translation of $A$. We also  denote by $|M|$ the number of pairwise nonisomorphic indecomposable summands of $M$, $\add M$ the subcategory consisting of direct summands of finite direct sums of $M$ for  $M\in \mod A$. Given an algebra $A=KQ/I$, let  $P_i$ be the indecomposable projective module,  $S_i$ the simple module, $e_i$  the primitive idempotent element  of an algebra corresponding to the point $i$.

\section{2. Preliminaries}
\subsection{Basic definitions}

In this section, we recall some basic definitions about ICE-closed subcategory in $\mod A$ and introduce the concept of epibrick in $\mod A$.
First of all, We give several conditions for a subcategory of $\mod A$, including closed under images, cokernels, extension, quotients and so on. 

\begin{definition}\label{def-2.1}
{\rm (\!\!\cite{2022ICE-closed})\ \
{\sl Let $A$ be an artin algebra and $\mathcal{T}$ a subcategory in $\mod A$.
	\begin{enumerate}
		\item[(1)]$\mathcal{T}$ is closed under images (resp. kernels, cokernels) if for every map $f: M\rightarrow N$ with $M, N \in \mathcal{T}$, we have Im $f\in \mathcal{T}$ (resp. Ker $f \in \mathcal{T}$, Coker $f \in \mathcal{T}$).
		\item[(2)]$\mathcal{T}$ is  closed under extensions if for every short exact sequence in $\mod A$
		$$0\rightarrow N_{1}\rightarrow N_{2}\rightarrow N_{3}\rightarrow 0$$
		with $N_{1}, N_{3} \in \mathcal{T}$, we have $N_{2} \in \mathcal{T}$ .
		\item[(3)]$\mathcal{T}$ is  closed under quotients if for every exact sequence in $\mod A$
		$$N_{1}\rightarrow N_{2}\rightarrow 0$$
		with $N_{1} \in \mathcal{T}$, we have $N_{2} \in \mathcal{T}$ .
	\end{enumerate}
}
}
\end{definition}

Then we can get the definitions of these subcategories.

\begin{definition}\label{def-2.2}
{\rm (\!\!\cite{2022ICE-closed})\ \
{\sl Let $A$ be an artin algebra and $\mathcal{T}$ a subcategory in $\mod A$.
	\begin{enumerate}
		\item[(1)]$\mathcal{T}$ is a torsion class if $\mathcal{T}$ is closed under quotients and extensions.
		\item[(2)]$\mathcal{T}$ is a wide subcategory if $\mathcal{T}$ is closed under kernels, cokernels, extensions.
		\item[(3)]$\mathcal{T}$ is an ICE-closed subcategory if $\mathcal{T}$ is closed under images, cokernels, extensions.
	\end{enumerate}
}
}
\end{definition}

\begin{corollary}\label{corollary-2.3}
	{\rm (\!\!\cite{2022ICE-closed})\ \
	{\sl  All torsion classes and wide subcategories are ICE-closed subcategories. 
	}
	}
\end{corollary}

{\bf Proof}\ \
	If $\mathcal{T}$ is a torsion class in $\mod A$, then $\mathcal{T}$ is closed under quotients and extensions.
	For every map $f: M\rightarrow N$ with $M, N \in \mathcal{T}$, We have 
	$$M\rightarrow Im f\rightarrow 0, N\rightarrow Coker f\rightarrow 0$$
	$Imf\in \mathcal{T}, Coker f\in \mathcal{T}$. Since $M, N \in \mathcal{T}$ and $\mathcal{T}$ is closed under extensions. $\mathcal{T}$ is an ICE-closed subcategory.
	
	In the same way, if $\mathcal{T}$ is a wide subcategory, then $\mathcal{T}$ is closed under kernels, cokernels and extensions.
	For every map $f: M\rightarrow N$ with $M, N \in \mathcal{T}$, We can get $Ker f\in \mathcal{T}$, $Coker f\in \mathcal{T}$.
	Then for map $g: N\rightarrow Coker f$ with $N, Coker f\in\mathcal{T}$, we have $Ker g\in\mathcal{T}$ and $ker g$ = $Im f$. That is $Im f \in\mathcal{T}$.
	$\mathcal{T}$ is an ICE-closed subcategory.
\\[8pt]

Next we give the definition of epibrick. 

\begin{definition}\label{def-2.4}{\rm (\!\!\cite{2021monobrick})\ \ {\sl Let $\mathcal{S}\in \mod A$.
			\begin{enumerate}
				\item[(1)]$\mathcal{S}$ is a brick  if $\End_{A}(\mathcal{S})$ is a division ring. The set of isoclasses of bricks in $\mod A$ is denoted by brick $A$.
				\item[(2)]A subset $\mathcal{S}\subseteq$ brick$A$ is called a semibrick if every morphism between elements of $\mathcal{S}$ is either zero or an isomorphism in $A$. The set of semibricks in $\mod A$ is denoted by sbrick $A$.
				\item[(3)]A subset $\mathcal{S}\subseteq$ brick$A$ is called a monobrick if every morphism between elements of $\mathcal{S}$ is either zero or an injection in $A$. The set of monobricks in $\mod A$ is denoted by mbrick $A$.
				\item[(4)]A subset $\mathcal{S}\subseteq$ brick$A$ is called an epibrick if every morphism between elements of $\mathcal{S}$ is either zero or a surjection in $A$. The set of epibricks in $\mod A$ is denoted by ebrick $A$.
			\end{enumerate}
	} }
\end{definition}

It is easy to know that every semibrick is a monobrick or epibrick.
By Schur's Lemma,every simple module is brick, and a set of isoclasses of simple modules is a semibrick.

Let $M\in \mod A$. The $one$-$point$ $extension$ of $A$ by $M_{A}$ is given by the following matrix algebra
$$B:=\left(\begin{array}{cc}A & 0 \\
	M_{A} & k\end{array}\right)$$
with the ordinary matrix and the multiplication induced by the module structure of $M_{A}$.
All $B$-modules can be seen as $(N_{A}, k^{n}, f)$, where $N_{A}\in \mod A$, $n \in \mathbb{N}$ and $f: k^{n}\otimes_k M_{A}\rightarrow N_{A}$.
The morphisms from $(N_{A}, k^{n_{1}},f_{1})$ to $(N'_{A}, k^{n_{2}}, f_{2})$ are pairs of ${(f, g)}$,where $f\in\Hom_A(N_{A},N'_{A})$ and $g\in\Hom(k^{n_{1}},k^{n_{2}})$,
such that the following diagram commutes,
$$\xymatrix{k^{n_{1}}\otimes_k M_{A}\ar[d]_{g\otimes M_{A}}\ar[rr]^{f_{1}}&&N_{A}\ar[d]^f\\
	k^{n_{2}}\otimes_k M_{A}\ar[rr]^{f_{2}}&&N'_{A}\\}$$

A sequence
$$0 \to {(N_{A}, k^{n_{1}}, f_{1})}\stackrel{{{ (h_1, g_1)}}}{\longrightarrow}
{(N'_{A}, k^{n_{2}}, f_{2})}\stackrel{{{(h_2, g_2)}}}{\longrightarrow}{(N''_{A}, k^{n_{3}}, f_{3})}\to 0$$
in $\mod B$ is exact if and only if 
$$0 \to N_{A}\stackrel{h_1}{\longrightarrow} N'_{A}\stackrel{h_2}{\longrightarrow} N''_{A} \to 0$$ is exact in $\mod A$ and
$$0 \to k^{n_{1}} \stackrel{g_1}{\longrightarrow} k^{n_{2}} \stackrel{g_2}{\longrightarrow} k^{n_{3}} \to 0$$ is exact in $\mod k$.

\section{3. Main Result}

In this section, we will give  ICE-closed subcategories (resp. epibricks) of $\mod B$ via an ICE-closed subcategory (resp. epibrick) of $\mod A$ in two different ways, where $B$ is one-point extension algebra of $A$ by an $A$-module $M_{A}$.

\begin{theorem}\label{3.1}
{\sl Let $B$ be the one-point extension algebra of $A$ by an $A$-module $M_{A}$ and $\mathcal{T}_{A}$ be an ICE-closed subcategory in $\mod A$.
	
	\begin{enumerate}
		\item[(1)]$\mathcal{T}_{B} := \{ (N_{A},0,0) | N_{A} \in \mathcal{T}_{A}\}$, $\mathcal{T}_{B}$ is an ICE-closed subcategory in $\mod B$.
		\item[(2)]$\mathcal{T}_{B} := \{ (N_{A},k^{n},f) ,(0,k^{n},0)| N_{A} \in \mathcal{T}_{A},n \in \mathbb{N},f: k^{n}\otimes_k M_{A}\rightarrow N_{A}\}$, $\mathcal{T}_{B}$ is an ICE-closed subcategory in $\mod B$.
\end{enumerate}}
\end{theorem}
{\bf Proof}\ \ 
	(1) Firstly, we check $\mathcal{T}_{B}$ is closed under extensions. Given an arbitrary short exact sequence in $\mod B$: 
	$0 \to {(N_{1}, 0, 0)}\longrightarrow{(N, k^{n}, f)}\longrightarrow{(N_{2}, 0, 0)}\to 0$, $(N_{1}, 0, 0)$, $(N_{2}, 0, 0)$ $\in \mathcal{T}_{B}$, we have 
	$0 \to N_{1}\longrightarrow N\longrightarrow N_{2} \to 0$ is exact in $\mod A$ and $0 \to 0 \longrightarrow k^{n} \longrightarrow 0 \to 0$ is exact in $\mod k$.
	Then $N\in \mathcal{T}_{A}$ and $n$ = 0. Since $N_{1}$, $N_{2}\in \mathcal{T}_{A}$ and $\mathcal{T}_{A}$ is closed under extension .
	Therefore, $(N, k^{n}, f)$ = $(N, 0, 0)$ $\in \mathcal{T}_{B}$.
	
	Secondly, we check $\mathcal{T}_{B}$ is closed under images. Given a map $F: (N_{1}, 0, 0)\rightarrow (N_{2}, 0, 0)$, $(N_{1}, 0, 0)$, $(N_{2}, 0, 0)\in \mathcal{T}_{B}$.
	$F$ = $(f,g)$, where $f: N_{1}\rightarrow N_{2}$, $g$ = 0.
	Obviously, $Im g$ = 0. $Im f\in \mathcal{T}_{A}$. Because $N_{1}$, $N_{2}\in \mathcal{T}_{A}$ and $\mathcal{T}_{A}$ is closed under images.
	$Im F$ = $(Im f, Im g, h)$ = $(Im f , 0 , 0) \in \mathcal{T}_{B}$.
	
	Finally, we check $\mathcal{T}_{B}$ is closed under cokernels. Given a map $F: (N_{1}, 0, 0)\rightarrow (N_{2}, 0, 0)$, $(N_{1}, 0, 0)$, $(N_{2}, 0, 0)\in \mathcal{T}_{B}$.
	$F$ = $(f,g)$, where $f: N_{1}\rightarrow N_{2}$, $g$ = 0.
	It is easy to know that $Coker g$ = 0 and $Coker f\in \mathcal{T}_{A}$. Because $N_{1}$, $N_{2}\in \mathcal{T}_{A}$ and $\mathcal{T}_{A}$ is closed under cokernels.
	$Coker F$ = $(Coker f, Coker g, h)$ = $(Coker f , 0 , 0) \in \mathcal{T}_{B}$.

	(2) Firstly, we check $\mathcal{T}_{B}$ is closed under extensions. Given an arbitrary short exact sequence in $\mod B$: 
	$0 \to {(N_{1}, k^{n_{1}}, f_{1})}\longrightarrow{(N_{2}, k^{n_{2}}, f_{2})}\longrightarrow{(N_{3}, k^{n_{3}}, f_{3})}\to 0$, $(N_{1},  k^{n_{1}}, f_{1})$, $(N_{3}, k^{n_{3}}, f_{3})$ $\in \mathcal{T}_{B}$, we have 
	$0 \to N_{1}\longrightarrow N_{2}\longrightarrow N_{3} \to 0$ is exact in $\mod A$ and $0 \to k^{n_{1}} \longrightarrow k^{n_{2}} \longrightarrow k^{n_{3}} \to 0$ is exact in $\mod k$.
	Then $N_{2} \in \mathcal{T}_{A}$ since $N_{1}$, $N_{3}\in \mathcal{T}_{A}$ and $\mathcal{T}_{A}$ is closed under extension. And $n_{2}$ = $n_{1}$ + $n_{3}$ $\in \mathbb{N}$.
	Therefore, $(N_{2}, k^{n_{2}}, f_{2})$ $\in \mathcal{T}_{A}$. In the same way, we can proof that $\{(0,k^{n},0)\}$ is closed under extensions.
	
	Secondly, we check $\mathcal{T}_{B}$ is closed under images. Given a map $F: (N_{1}, k^{n_{1}}, f_{1})\rightarrow (N_{2}, k^{n_{2}}, f_{2})$, $(N_{1}, k^{n_{1}}, f_{1})$, $(N_{2}, k^{n_{2}}, f_{2})\in \mathcal{T}_{B}$.
	$F$ = $(f,g)$, where $f: N_{1}\rightarrow N_{2}$, $g$ = $k^{n_{1}} \rightarrow k^{n_{2}}$.
	$N_{1}$, $N_{2}\in \mathcal{T}_{A}$ and $\mathcal{T}_{A}$ is closed under images. So $Im f \in \mathcal{T}_{A}$.
	$Im g$ is subspace of $k^{n_{2}}$. Then $Im g$ is $n$ dimensional vector space, $n \in \mathbb{N}$. 
	$Im F$ = $(Im f, Im g, h)$ $\in \mathcal{T}_{B}$, $h:Im g\otimes_k M_{A}\rightarrow Im f$.
	Similarly, we can proof that $\{(0,k^{n},0)\}$ is closed under images($f$ = 0).
	
	Finally, we check $\mathcal{T}_{B}$ is closed under cokernels. Given a map $F: (N_{1}, k^{n_{1}}, f_{1})\rightarrow (N_{2}, k^{n_{2}}, f_{2})$, $(N_{1}, k^{n_{1}}, f_{1})$, $(N_{2}, k^{n_{2}}, f_{2})\in \mathcal{T}_{B}$.
	$F$ = $(f,g)$, where $f: N_{1}\rightarrow N_{2}$, $g$ = $k^{n_{1}} \rightarrow k^{n_{2}}$.
	$N_{1}$, $N_{2}\in \mathcal{T}_{A}$ and $\mathcal{T}_{A}$ is closed under cokernels. So $Coker f \in \mathcal{T}_{A}$.
	Obviously $Coker g$ is $n$ dimensional vector space, $n \in \mathbb{N}$.
	$Coker F$ = $(Coker f, Coker g, h) \in \mathcal{T}_{B}$, $h:Coker g\otimes_k M_{A}\rightarrow Coker f$.
	Similarly, we can proof that $\{(0,k^{n},0)\}$ is closed under cokernels($f$ = 0).
\\[8pt]

\begin{example}{\rm
	\begin{enumerate}
		\item[(1)]$B$ := $KQ_{B}$, $Q_{B}$ : $1\stackrel{\alpha}{\rightarrow} 2\stackrel{\beta}{\rightarrow} 3$.Let $A$ := $KQ_{A}$, $Q_{A}$ : $2\stackrel{\beta}{\rightarrow} 3$, $M_{A}$ = $\left\langle\alpha\right\rangle$ $\cong$ $P_{2}$. Then $B:=\left(\begin{array}{cc}A & 0 \\
			M_{A} & k\end{array}\right)$.	
		The irreducible representations of $A$ are $P_{2}$ : $k\rightarrow k$, $S_{2}$ : $k\rightarrow 0$, $S_{3}$ : $0\rightarrow k$. The ICE-closed subcategories in $\mod A$ are 
		add$\left\{\begin{matrix}	2 \\ 3	\end{matrix},2,3\right\}$, add$\left\{\begin{matrix}	2 \\ 3	\end{matrix},2\right\}$, add$\left\{2\right\}$, add$\left\{3\right\}$, add$\left\{0\right\}$.
		Then we can get ICE-closed subcategories in $\mod B$ :
		add$\left\{\begin{matrix}	2 \\ 3	\end{matrix},2,3\right\}$, add$\left\{\begin{matrix}	2 \\ 3	\end{matrix},2\right\}$, add$\left\{2\right\}$, add$\left\{3\right\}$, add$\left\{0\right\}$,
		add$\left\{\begin{matrix}	1\\ 2 \\ 3	\end{matrix},1,\begin{matrix}	2 \\ 3	\end{matrix},\begin{matrix}	1\\ 2\end{matrix},2,3\right\}$, add$\left\{\begin{matrix}	1\\ 2 \\ 3	\end{matrix},1,\begin{matrix}	2 \\ 3	\end{matrix},\begin{matrix}	1\\ 2\end{matrix},2\right\}$, add$\left\{ \begin{matrix}	1\\ 2\end{matrix},1,2\right\}$, add$\{1,3\}$, add$\{1\}$ by Theorem 3.1.
		\item[(2)]$B$ := $KQ_{B}$, $Q_{B}$ : $1\stackrel{\alpha}{\rightarrow} 2\stackrel{\beta}{\rightarrow} 3$ with relation $\alpha\beta = 0$.
		Let $A$ := $KQ_{A}$, $Q_{A}$ : $2\stackrel{\beta}{\rightarrow} 3$, $M_{A}$ = $\left\langle\alpha\right\rangle$ $\cong$ $S_{2}$. Then $B:=\left(\begin{array}{cc}A & 0 \\
			M_{A} & k\end{array}\right)$.
		The irreducible representations of $A$ are the same as $(1)$.
		The ICE-closed subcategories in $\mod A$ are also identical to $(1)$. However, the ICE-closed subcategories in $\mod B$ are add$\left\{\begin{matrix}	2 \\ 3	\end{matrix},2,3\right\}$, add$\left\{\begin{matrix}	2 \\ 3	\end{matrix},2\right\}$, add$\left\{2\right\}$, add$\left\{3\right\}$, add$\left\{0\right\}$,
		add$\left\{\begin{matrix}	2 \\ 3	\end{matrix},1,\begin{matrix}	1\\ 2\end{matrix},2,3\right\}$, add$\left\{ \begin{matrix}	2 \\ 3	\end{matrix},1,\begin{matrix}	1\\ 2\end{matrix},2\right\}$, add$\left\{ \begin{matrix}	1\\ 2\end{matrix},1,2\right\}$, add$\{1,3\}$, add$\{1\}$ by Theorem 3.1.
	\end{enumerate}
}
\end{example}

\begin{remark}\label{3.3} {\rm
	Applying Theorem 3.1, we can give a part of ICE-closed subcategories in $\mod B$. But more computation is required to give all the ICE-closed subcategories in $\mod B$.}
\end{remark}

\begin{corollary} {\sl
	The number of ICE-closed subcategories in $\mod B$ and the number of ICE-closed subcategories in $\mod A$ are denoted respectively as $m$, $n$.
	Then we have : $$m \geq 2n.$$}
\end{corollary}

\begin{theorem} {\sl
	Let $B$ be the one-point extension algebra of $A$ by an $A$-module $M_{A}$ and $\mathcal{S}_{A}$ be an epibrick in $\mod A$.
	
	\begin{enumerate}
		\item[(1)]$\mathcal{S}_{B} := \{ (s,0,0) | s \in \mathcal{S}_{A}\}$, $\mathcal{S}_{B}$ is an epibrick in $\mod B$.
		\item[(2)]$\mathcal{S'}_{B} := \{ (s,0,0) ,(0,k,0)| s \in \mathcal{S}_{A}\}$, $\mathcal{S'}_{B}$ is an epibrick in $\mod B$.
	\end{enumerate}
}
\end{theorem}
{\bf Proof}\ \ 
		(1)For an arbitrary morphism $F : w_{1}\rightarrow w_{2}$, where $w_{1}=(s_{1},0,0)$, $w_{2}=(s_{2},0,0) \in \mathcal{S}_{B}$, it is easy to know $F=(f,0)$ with $f:s_{1}\rightarrow s_{2}$. So $F \cong f$. $f$ is either zero or a surjection. Since $s_{1}$, $s_{2} \in \mathcal{S}_{A}$ and $\mathcal{S}_{A}$ is an epibrick in $\mod A$.
		Therefore $F$ is either zero or a surjection. $\mathcal{S}_{B}$ is an epibrick in $\mod B$.
		
		(2)For an arbitrary morphism $F : w_{1}\rightarrow w_{2}$, where $w_{1}=(s,0,0)$, $w_{2}=(0,k,0) \in \mathcal{S}_{B}$, $F=(f,g)$ with $f:s\rightarrow 0$, $g:0\rightarrow k$. So $F=0$. According to $(1)$, $\mathcal{S}_{B}$ is an epibrick in $\mod B$ and ${(0,k,0)}$ is also an epibrick in $\mod B$. Therefore, $\mathcal{S'}_{B}$ is an epibrick in $\mod B$.
\\[8pt]
\begin{remark}\label{3.6} {\rm
	Let $B$ be the one-point extension algebra of $A$ by an $A$-module $M_{A}$ and $\mathcal{S}_{A}$ be a monobrick in $\mod A$.
	
	\begin{enumerate}
		\item[(1)]$\mathcal{S}_{B} := \{ (s,0,0) | s \in \mathcal{S}_{A}\}$, $\mathcal{S}_{B}$ is an monobrick in $\mod B$
		\item[(2)]$\mathcal{S'}_{B} := \{ (s,0,0) ,(0,k,0)| s \in \mathcal{S}_{A}\}$, $\mathcal{S'}_{B}$ is an monobrick in $\mod B$
	\end{enumerate}
}
\end{remark}
\begin{example} {\rm
	$B$ := $KQ_{B}$, $Q_{B}$ : 
	$$	\begin{tikzpicture}[thick]
		\node (1) at (0,0) {$1$};
		\node (2) at (1.5,0) {$2$};
		\node (4) at (-15:3) {$4$};
		\node (3) at (15:3) {$3$};
		\draw (0.75,0.05) node[above] {$\alpha$};
		\draw (10:2.2) node[above] {$\beta$};
		\draw (-10:2.2) node[above] {$\gamma$};
		\draw[->] (1)--(2);
		\draw[->] (2)--(4);
		\draw[<-] (2)--(3);
	\end{tikzpicture}$$
	Let $A$ := $KQ_{A}$, $Q_{A}$ : $4\stackrel{\gamma}{\rightarrow}2\stackrel{\beta}{\rightarrow} 3$, $M_{A}$ = $\left\langle\alpha\right\rangle$ = $k\{\alpha,\alpha\beta\}$ $\cong$ $P_{2}$ : $0\rightarrow k \rightarrow k$.Then $B:=\left(\begin{array}{cc}A & 0 \\
		M_{A} & k\end{array}\right)$.
	The irreducible representations of $A$ : $4$, $2$, $3$, $\begin{matrix}	4\\ 2\end{matrix}$, $\begin{matrix}	2\\ 3\end{matrix}$, $\begin{matrix}	4\\ 2\\3 \end{matrix}$.
	The epibricks in $\mod A$ are:
	$\left\{ 4 \right\}$, $\left\{ 4,2 \right\}$, $\left\{ 4,3 \right\}$, $\left\{ 4,\begin{matrix}	4\\ 2\end{matrix} \right\}$, $\left\{ 4,\begin{matrix}	2\\ 3\end{matrix} \right\}$, $\left\{ 4,\begin{matrix}4\\2\\ 3\end{matrix} \right\}$, $\left\{ 4,2,3 \right\}$, $\left\{ 4,2,\begin{matrix}	2\\ 3\end{matrix} \right\}$, $\left\{ 4,2,\begin{matrix}	4\\2\\ 3\end{matrix} \right\}$, $\left\{ 4,3,\begin{matrix}	4\\ 2\end{matrix} \right\}$, $\left\{ 4,\begin{matrix}	4\\ 2\end{matrix},\begin{matrix}	4\\ 2\\3\end{matrix} \right\}$, $\left\{ 2 \right\}$, 	$\left\{ 2,3 \right\}$, $\left\{ 2,\begin{matrix}	2\\ 3\end{matrix} \right\}$, $\left\{ 2,\begin{matrix}	4\\ 2\\3\end{matrix} \right\}$, $\left\{ 2,\begin{matrix}	4\\ 2\\3\end{matrix},\begin{matrix}	4\\ 2\end{matrix} \right\}$, $\{3\}$, $\left\{ 3,\begin{matrix}	4\\ 2\end{matrix} \right\}$, $\left\{\begin{matrix}	4\\ 2\end{matrix},\begin{matrix}	4\\ 2\\3\end{matrix} \right\}$, $\left\{\begin{matrix}	2\\3\end{matrix}\right\}$, $\left\{\begin{matrix}	4\\ 2\\3\end{matrix}\right\}$, $\{0\}$.
	
	Then we can get some epibricks in $mod B$ by Theorem 3.5:
	
	$\left\{ 4 \right\}$, $\left\{ 4,2 \right\}$, $\left\{ 4,3 \right\}$, $\left\{ 4,\begin{matrix}	4\\ 2\end{matrix} \right\}$, $\left\{ 4,\begin{matrix}	2\\ 3\end{matrix} \right\}$, $\left\{ 4,\begin{matrix}4\\2\\ 3\end{matrix} \right\}$, $\left\{ 4,2,3 \right\}$, $\left\{ 4,2,\begin{matrix}	2\\ 3\end{matrix} \right\}$, $\left\{ 4,2,\begin{matrix}	4\\2\\ 3\end{matrix} \right\}$, $\left\{ 4,3,\begin{matrix}	4\\ 2\end{matrix} \right\}$, $\left\{ 4,\begin{matrix}	4\\ 2\end{matrix},\begin{matrix}	4\\ 2\\3\end{matrix} \right\}$, $\left\{ 2 \right\}$, 	$\left\{ 2,3 \right\}$, $\left\{ 2,\begin{matrix}	2\\ 3\end{matrix} \right\}$, $\left\{ 2,\begin{matrix}	4\\ 2\\3\end{matrix} \right\}$, $\left\{ 2,\begin{matrix}	4\\ 2\\3\end{matrix},\begin{matrix}	4\\ 2\end{matrix} \right\}$, $\{3\}$, $\left\{ 3,\begin{matrix}	4\\ 2\end{matrix} \right\}$, $\left\{\begin{matrix}	4\\ 2\end{matrix},\begin{matrix}	4\\ 2\\3\end{matrix} \right\}$, $\left\{\begin{matrix}	2\\3\end{matrix}\right\}$, $\left\{\begin{matrix}	4\\ 2\\3\end{matrix}\right\}$, $\{0\}$,
	$\left\{ 4,1 \right\}$, $\left\{ 4,2,1 \right\}$, $\left\{ 4,3,1 \right\}$, $\left\{ 4,\begin{matrix}	4\\ 2\end{matrix},1 \right\}$, $\left\{ 4,\begin{matrix}	2\\ 3\end{matrix},1 \right\}$, $\left\{ 4,\begin{matrix}4\\2\\ 3\end{matrix},1 \right\}$, $\left\{ 4,2,3,1 \right\}$, $\left\{ 4,2,\begin{matrix}	2\\ 3\end{matrix},1 \right\}$, $\left\{ 4,2,\begin{matrix}	4\\2\\ 3\end{matrix},1 \right\}$, $\left\{ 4,3,\begin{matrix}	4\\ 2\end{matrix},1 \right\}$, $\left\{ 4,\begin{matrix}	4\\ 2\end{matrix},\begin{matrix}	4\\ 2\\3\end{matrix},1 \right\}$, $\left\{ 2,1 \right\}$, 	$\left\{ 2,3,1 \right\}$, $\left\{ 2,\begin{matrix}	2\\ 3\end{matrix},1 \right\}$, $\left\{ 2,\begin{matrix}	4\\ 2\\3\end{matrix},1 \right\}$, $\left\{ 2,\begin{matrix}	4\\ 2\\3\end{matrix},\begin{matrix}	4\\ 2\end{matrix},1 \right\}$, $\{3,1\}$, $\left\{ 3,\begin{matrix}	4\\ 2\end{matrix},1 \right\}$, $\left\{\begin{matrix}	4\\ 2\end{matrix},\begin{matrix}	4\\ 2\\3\end{matrix},1 \right\}$, $\left\{\begin{matrix}	2\\3\end{matrix},1\right\}$, $\left\{\begin{matrix}	4\\ 2\\3\end{matrix},1\right\}$, $\{1\}$.
}
\end{example}
\begin{remark}\label{3.8} {\rm
	
	Applying Theorem 3.5, we can give a part of epibricks in $\mod B$. But more computation is required to give all the epibricks in $\mod B$.
}
\end{remark}

\begin{corollary} {\sl
	The number of epibricks in $\mod B$ and the number of epibricks in $\mod A$ are denoted respectively as $m$, $n$.
	Then we have : $$m \geq 2n.$$
}
\end{corollary}
{\bf Applications.} 
Let $\Lambda$ be an algebra and $M\in\mod\Lambda$. $M$ is $\tau$-tilting if $\Hom_\Lambda(M,\tau M)=0$ and  $|M|=|\Lambda|$. $M$ is  support $\tau$-tilting if it is a $\tau$-tilting $\Lambda/\Lambda e\Lambda$-module for some idempotent $e$ of $\Lambda$. Enomoto shown that every  every  functorially finite wide subcategory  $\mathcal{W}$ is  equivalent to a module category (i.e, there is an algebra $\Gamma$ such that   $\mathcal{W}$ is equivalent to $\mod\Gamma$), and then he introduced the definition of wide $\tau$-tilting modules as follows.

\begin{definition}{\rm (\cite{2022ICE-closed})}{\sl
	\begin{enumerate}
		\item[(1)] Given a  functorially finite wide subcategory  $\mathcal{W}$ of $\mod\Lambda$ and $M\in\mathcal{W}$, fix a equivalent $F:\mathcal{W}\simeq\mod\Gamma$. We say $M$ is $\tau_{\mathcal{W}}$-tilting if $F(M)$ is a $\tau$-tilting $\Gamma$-module.
		\item[(2)]A $\Lambda$-module $M$ is called {\it wide~$\tau$-tilting}  if there is a   functorially finite wide subcategory  $\mathcal{W}$ of $\mod\Lambda$ such that $M$ is  $\tau_{\mathcal{W}}$-tilting. The set of all wide $\tau$-tilting $\Lambda$-modules will be denoted by {\text w}$\tau$-tilt $\Lambda$.
	\end{enumerate}
}
\end{definition}

Suppose that $A$, $B$ are Nakayama algebras and $B$ is the one-point extension of $A$ by an $A$-module $M_{A}$. In \cite{Gao}, the authors get the following bijections:
$$\begin{tikzpicture}[thin]
	\node (1) at (0,0) {w$\tau$$-$tilt $\Lambda$};
	\node (2) at (3,0) {ice $\Lambda$};
	\node (3) at (6,0) {ebrick $\Lambda$};
	
	\draw (1.75,0.01) node[above] {$\tiny cok(-)$};
	\draw (1.75,-0.01) node[below] {$\tiny P(-)$};
	\draw (4.4,0.01) node[above] {$\tiny Sim(-)$};
	\draw (4.4,-0.01) node[below] {$\tiny Filt(-)$};
	\draw[->] (1,0.05)--(2.5,0.05);
	\draw[<-] (1,-0.05)--(2.5,-0.05);
	\draw[->] (3.5,0.05)--(5.25,0.05);
	\draw[<-] (3.5,-0.05)--(5.25,-0.05);
\end{tikzpicture}$$
where $\Lambda$ is either $A$ or $B$, cok($M$) denote the subcategory of mod $A$ consisting of cokernels of morphisms in add $M$, Filt($\mathcal{S}$) denote the minimal Extension-closed subcategory which contains $\mathcal{S}$ for $\mathcal{S} \in$ ebrick $A$, Sim($B$) denote the set of all simple object of ice $B$, P($\mathcal{C}$) denote the maximal Ext-projective object of $\mathcal{C}$.
Then we have two different ways to construct wide $\tau$-tilting $B$-modules from wide $\tau$-tilting $A$-modules as follows:
$$\begin{tikzpicture}[thin]
	\node (2) at (0,2) {ice $A$};
	\node (3) at (0,4) {w$\tau$-tilt $A$};
	\node (5) at (4,2) {ice $B$};
	\node (6) at (4,4) {w$\tau$-tilt $B$};
	
	\draw (2,2.05) node[above] {Theorem 3.1};
	\draw (0,3) node[left] {$cok( - )$};
	\draw (4,3) node[right] {$P( - )$};
	\draw[->] (3)--(2);
	\draw[->] (5)--(6);
	\draw[->] (2)--(5);
\end{tikzpicture}$$
and
$$\begin{tikzpicture}[thin]
	\node (2) at (0,2) {ebrick $A$};
	\node (3) at (0,4) {w$\tau$-tilt $A$};
	\node (5) at (4,2) {ebrick $B$};
	\node (6) at (4,4) {w$\tau$-tilt $B$};
	
	\draw (2,2.05) node[above] {Theorem 3.5};
	\draw (0,3) node[left] {${\tiny Sim(cok(-))}$};
	\draw (4,3) node[right] {${\tiny P(Filt(-))}$};
	\draw[->] (3)--(2);
	\draw[->] (5)--(6);
	\draw[->] (2)--(5);
\end{tikzpicture}$$

\begin{example} {\rm
	$B$ := $KQ_{B}$, $Q_{B}$ : $1\stackrel{\alpha}{\rightarrow} 2\stackrel{\beta}{\rightarrow} 3$.
	Let $A$ := $KQ_{A}$, $Q_{A}$ : $2\stackrel{\beta}{\rightarrow} 3$, $M_{A}$ = $\left\langle\alpha\right\rangle$ $\cong$ $P_{2}$. Then $B:=\left(\begin{array}{cc}A & 0 \\
		M_{A} & k\end{array}\right)$.
	\begin{enumerate}
		\item[1.]\begin{enumerate}
			\item [(1)]w$\tau$-tilt $A\subseteq$ w$\tau$-tilt $B$ by Theorem 3.1(1).
			\item [(2)]We list w$\tau$-tilt $A$, ice $A$, ice $B$ and w$\tau$-tilt $B$ in table 1 by Theorem 3.1(2).
			\newpage
			\begin{center}
				\begin{tabular}{|c|c|c|c|}
					\hline
					w$\tau$-tilt $A$ & ice $A$ & ice $B$ & w$\tau$-tilt $B$ \\
					\hline
					0&add\{0\}&add$\{1\}$&1\\
					\hline
					2&add\{2\}&add$\left\{\begin{matrix}1\\2\end{matrix}, 1, 2\right\}$&2 $\oplus$ $\begin{matrix}1\\2\end{matrix}$ \\
					\hline
					3&add\{3\}&add\{1, 3\}&1 $\oplus$ 3 \\
					\hline
					$\begin{matrix}2\\3\end{matrix}$&add$\left\{\begin{matrix}2\\3\end{matrix}\right\}$&add$\left\{\begin{matrix}1\\2\\3\end{matrix}, 1, \begin{matrix}2\\3\end{matrix}\right\}$&$\begin{matrix}2\\3\end{matrix}$ $\oplus$ $\begin{matrix}1\\2\\3\end{matrix}$ \\
					\hline
					$\begin{matrix}2\\3\end{matrix}$ $\oplus$ 2&add$\left\{\begin{matrix}2\\3\end{matrix},2\right\}$&add$\left\{\begin{matrix}1\\2\\3\end{matrix},1,\begin{matrix}2\\3\end{matrix},\begin{matrix}1\\2\end{matrix},2\right\}$&$\begin{matrix}1\\2\\3\end{matrix}$ $\oplus$ $\begin{matrix}2\\3\end{matrix}$ $\oplus$ 2 \\
					\hline
					$\begin{matrix}2\\3\end{matrix}$ $\oplus$ 3&add$\left\{\begin{matrix}2\\3\end{matrix},3,2\right\}$&$\mod$ $B$&$\begin{matrix}1\\2\\3\end{matrix}$ $\oplus$ $\begin{matrix}2\\3\end{matrix}$ $\oplus$ 3 \\
					\hline
				\end{tabular}
		\end{center}
		\begin{center}
			{\small Table 1~~w$\tau$-tilt $A$ ice $A$ ice $B$ w$\tau$-tilt $B$}
		\end{center}
		\end{enumerate}
		\item[2.]\begin{enumerate}
			\item [(1)]w$\tau$-tilt $A\subseteq$ w$\tau$-tilt $B$ by Theorem 3.5(1).
			\item [(2)]We list w$\tau$-tilt $A$, ebrick $A$, ebrick $B$ and w$\tau$-tilt $B$ in table 2 by Theorem 3.5(2).
             \begin{center}
				\begin{tabular}{|c|c|c|c|}
					\hline
					w$\tau$-tilt $A$ & ebrick $A$ & ebrick $B$ & w$\tau$-tilt $B$ \\
					\hline
					0&\{0\}&\{1\}&1\\
					\hline
					2&\{2\}&\{1,2\}&2 $\oplus$ $\begin{matrix}1\\2\end{matrix}$ \\
					\hline
					3&\{3\}&\{1, 3\}&1 $\oplus$ 3 \\
					\hline
					$\begin{matrix}2\\3\end{matrix}$&$\left\{\begin{matrix}2\\3\end{matrix}\right\}$&$\left\{1, \begin{matrix}2\\3\end{matrix}\right\}$&$\begin{matrix}2\\3\end{matrix}$ $\oplus$ $\begin{matrix}1\\2\\3\end{matrix}$ \\
					\hline
					$\begin{matrix}2\\3\end{matrix}$ $\oplus$ 2&$\left\{\begin{matrix}2\\3\end{matrix},2\right\}$&$\left\{1,\begin{matrix}2\\3\end{matrix},2\right\}$&$\begin{matrix}1\\2\\3\end{matrix}$ $\oplus$ $\begin{matrix}2\\3\end{matrix}$ $\oplus$ 2 \\
					\hline
					$\begin{matrix}2\\3\end{matrix}$ $\oplus$ 3&\{2,3\}&\{1,2,3\}&$\begin{matrix}1\\2\\3\end{matrix}$ $\oplus$ $\begin{matrix}2\\3\end{matrix}$ $\oplus$ 3 \\
					\hline
				\end{tabular}
			\end{center}
			\begin{center}
				{\small Table 2~~w$\tau$-tilt $A$ ebrick $A$ ebrick $B$ w$\tau$-tilt $B$}
			\end{center}
		\end{enumerate}
	\end{enumerate}
}
\end{example}

{\bf Acknowledgements.} The authors would like to thank Yanhong Bao for helpful discussions. The second author was supported by NSFC (Grant No.12301041 ). 

\vspace{0.6cm}
%


%

{\footnotesize
\begin{thebibliography}{99}
\addtolength{\itemsep}{-0.6ex}

\bibitem{2020semibrick} S. ASAI. {\sl Semibricks}. Int.Math.Res.Not.IMRN., 2020, {\bf 16}: 4993--5054.

\bibitem{1976re} C. M. RINGEL. {\sl Representations of $K$-species and bimodules}. J.Algebra., 1976, {\bf 41}(2): 269--302.

\bibitem{2022ICE-closed} H. ENOMOTO. {\sl Rigid modules and ICE-closed subcategories in quiver representations}. Journal of Algebra., 2022, {\bf 594}: 364--388.

\bibitem{2021monobrick} H. ENOMOTO. {\sl Monobrick, a uniform approach to torsion-free classes and wide subcategories}. Advances in Mathematics., 2021, {\bf 393}: 108--113.

\bibitem{Gao} Hanpeng GAO, Dajun LIU. {\sl A note on wide $\tau$-tilting modules and epibricks}. prepublication.

\bibitem{1962france} P. GABRIEL. {\sl Des cat\'{e}gories ab\'{e}liennes}. Bull.Soc.Math.France., 1962, {\bf 90}: 323--448.

\bibitem{2014tiltingtheory} T. ADACHI, O. IYAMA, I. REITEN. {\sl $\tau$-tilting theory}. Compos.Math., 2014, {\bf 150}(3): 415--452.




\end{thebibliography}

}

\end{document}